\documentstyle[amssymb,12pt]{amsart}

\ifx\mathrm\undefined\let\mathrm\bf\fi
\ifx\mathbf\undefined\let\mathbf\bf\fi

\let\leq\leqslant 
\let\tsize\textstyle \def
\Sum{\sum\limits}

\let\al\alpha

 \let\eps\varepsilon \let\epsilon\eps

\let\la\lambda \let\La\Lambda
\let\om\omega \let\Om\Omega
 \let\phi\varphi

\let\tht\theta

\newcommand{\Zb}{{{\mathcal Z}}}

\newcommand{\half}{\frac12}
\newcommand{\Z}{{\Bbb Z}}

\newcommand{\C}{{\Bbb C}}
   
\newcommand{\Ref}[1]{{$($\ref{#1}$)$}}
\newcommand{\bean}{\begin{eqnarray}}
\newcommand{\eean}{\end{eqnarray}}
\newcommand{\be}{\begin{displaymath}}
\newcommand{\ee}{\end{displaymath}}
\newcommand{\bea}{\begin{eqnarray*}}   
\newcommand{\eea}{\end{eqnarray*}}

\newcommand{\h}{{{\frak h\,}}}
\newcommand{\Id}{{\operatorname{Id}}}

\newcommand{\res}{{\operatorname{res}}}
\renewcommand{\Im}{{\operatorname{Im}}\,}
\renewcommand{\Re}{{\operatorname{Re}}\,}

\newcommand{\T}{\otimes}

\newcommand{\vs}{\vspace{1.5\baselineskip}}

\newtheorem
{thm}{Theorem}

\newtheorem
{lemma}[thm]{Lemma}

\newcommand{\th}{\theta}
\newcommand{\End}{{\operatorname{End}}}

\begin{document}
\title[The qKZB Equation in Tensor Products of Finite Dimensional Modules]
{Solutions of the qKZB Equation in Tensor Products of finite dimensional modules over the elliptic quantum group $E_{\tau,\eta}sl_2$}
\author{E. Mukhin and A. Varchenko}
\maketitle
\vskip-.5\baselineskip
\centerline{${}^*${\it Department of Mathematics,
University of North Carolina at Chapel Hill,}}
\centerline{\it Chapel Hill, NC 27599-3250, USA}
\centerline{{\it E-mail addresses:} {\rm mukhin@@math.unc.edu,
av@@math.unc.edu}}
\medskip

\centerline{December, 1997}

\begin{abstract}

We consider the quantized Knizhnik-Zamolodchikov-Bernard
difference equation (qKZB) with step $p$ and values in a tensor product of finite dimensional evaluation modules over the elliptic quantum group $E_{\tau,\eta}(sl_2)$, the equation defined in terms of elliptic dynamical R-matrices.
We solve the equation in terms of multidimensional q-hypergeometric
integrals and describe its monodromy properties. We identify  the space of solutions of the qKZ equation with the space of functions with values in the tensor product of the corresponding modules over the elliptic quantum group $E_{p,\eta}(sl_2)$.
\end{abstract}

\section{Introduction}

In this paper we solve the system of elliptic quantum
Knizhnik--Zamolodchikov-Bernard (qKZB) difference equations associated with
the elliptic quantum group $E_{\tau,\eta}(sl_2)$ with values in finite dimensional evaluation representations and describe the monodromy
properties of solutions.

The qKZB equation \cite{F} is a quantum deformation of the KZB differential
equation obeyed by correlation functions of the Wess-Zumino-Witten model on
tori. The qKZB equation has the form
\be
\Psi(z_1,\dots,z_j+p,\dots,z_n)=
K_j(\la,\tau,p,\eta;z_1,\dots,z_n;\eta\La_1,\dots,\eta\La_n)\Psi(z_1,\dots,z_n).
\ee
The unknown function $\Psi$ takes values in a space of vector valued functions
of a complex variable $\lambda$, and the operators $K_j$ are
expressed in terms of $R$-matrices of the elliptic quantum group
$E_{\tau,\eta}(sl_2)$. The parameters of this system of equations are $\tau$ (the period
of the elliptic curve), $\eta$ (``Planck's constant''), $p$ (the step) and $n$
``highest weights'' $\Lambda_1,\dots,\Lambda_n\in\C$. 

In the trigonometric limit $\tau\to i\infty$, the qKZB equation reduces to
the trigonometric qKZ equation \cite{FR} obeyed by correlation functions of
statistical models and form factors of integrable quantum field theories in
1+1 dimensions \cite{JM, S}.

The KZB differential equation can be obtained in the semiclassical limit: $\eta\to 0$,
$p\to 0$, $p/\eta$ finite.

Solutions of the qKZB equation with values in a tensor product 
$V_{\La_1}\T\ldots\T V_{\La_n}$
of 
evaluation Verma modules over the elliptic quantum group $E_{\tau,\eta}(sl_2)$ with generic highest weights 
$\La_1,\dots,\La_n$
were constructed in \cite{FTV1}, \cite{FTV}. The solutions have the form
\bea
\Psi(z)\,=\, \sum_{k_1,...,k_n} I_{k_1,...,k_n}(z) f^{k_1}v_1 \otimes
... \otimes  f^{k_n}v_n,
\eea
where $\{f^{k_1}v_1\T\ldots\T f^{k_n}v_n\}\in V_{\La_1}\T ...\T V_{\La_n}$ is 
the standard basis in the tensor product of the Verma modules,
and the coefficients $I_{k_1,...,k_n}(z)$ are given by suitable multidimensional hypergeometric integrals. 

If weights become positive integers, i.e. if the evaluation Verma modules of the tensor product $V_{\La_1}\T ...\T V_{\La_n}$
become reducible , then some of the coefficients in the sum become divergent.

In this paper we prove that the restricted sum
\bea
\Psi^0(z)\,=\, 
\sum_{k_1\le\La_1 ,...,k_n\le\La_n} I_{k_1,...,k_n}(z) 
f^{k_1}v_1 \T\ldots\T f^{k_n}v_n
\eea
remains well defined. Moreover, we show that the sum $\Psi^0(z)$
defines a solution of the qKZB equation with values
in the tensor product $L_{\la_1}\T ... \T L_{\la_n}$ of finite dimensional evaluation modules over the elliptic group $E_{\tau,\eta}(sl_2)$. We use the method developed in \cite{MV} to construct solutions of the rational qKZ equation with values in irreducible $sl_2$ modules.

This result allows us to obtain a description of the monodromy properties of
solutions in a way parallel to \cite{FTV}. The monodromy transformations of solutions are described in terms of $R$-matrices associated
with pairs of representations of the "dual" elliptic quantum group
$E_{p,\eta}(sl_2)$, where $p$ is the step of the difference equation.
The description of the monodromy is analogous to the Kohno-Drinfeld description
\cite{K}, \cite{D} of the monodromy group of solutions of the KZ differential
equations associated to a simple Lie algebra in terms of the corresponding
quantum group.

We show that the residues of divergent hypergeometric solutions of the qKZB equation are again hypergeometric solutions of the qKZB equation with new weights up to some explicit scalar factors.

The paper is organized as follows.

In Section \ref{weight fun} we recall a construction of the space of elliptic weight functions. In Section \ref{QKZB} we describe the elliptic $R$-matrices and the qKZB equation. In Section \ref{integrals} we consider hypergeometric integrals. This Section contains the main results of this paper, Theorems \ref{an.cont} and \ref{functor}. In Section \ref{solutions} we use the results from Section \ref{integrals} and \cite{FTV} to describe solutions of the qKZB equation and their monodromy properties.

\section{Elliptic weight functions}\label{weight fun}

\subsection{The phase function}

Fix natural numbers $n,l$. Fix complex parameters $\tau$, $\eta$, $p$ such that $\Im\tau>0$, $\Im\eta<0$ and $\Im p>0$. Fix $\La=(\Lambda_1,\dots, \Lambda_n),\, z=(z_1,\dots, z_n)\in\C^n$. Set $r=e^{2\pi ip}$, $q=e^{2\pi i\tau}$, $a_j=\eta\La_j,\, j=1,\dots,n$ . 

Recall that the Jacobi theta function
\be
\tht_\tau(t)=-\sum_{j\in\Z}
e^{\pi i(j+\half)^2\tau+2\pi i(j+\half)(t+\half)},
\ee
has multipliers $-1$ and $-\exp(-2\pi it-\pi i\tau)$ as $t\to t+1$ and
$t\to t+\tau$, respectively. It is an odd entire function whose zeros are
simple and lie on the lattice $\Z+\tau \Z$.

Define one variable \emph{phase function} $\Omega(t,a)$ by the convergent infinite product

\bean\label{phase}
\Omega(t,a):=\Omega(t,\tau,p;a)=\prod_{j=0}^{\infty}\prod_{k=0}^{\infty}
\frac{(1-r^jq^ke^{2\pi i(t-a)})(1-r^{j+1}q^{k+1}e^{-2\pi i(t+a)})}
{(1-r^jq^ke^{2\pi i(t+a)})(1-r^{j+1}q^{k+1}e^{-2\pi i(t-a)})}\,.
\eean

The one variable phase function $\Omega(t,a)$ has the properties
\be\label{ph}
\Omega(t+p,a)=r^a\,\frac{\theta_\tau(t+a)}{\theta_\tau(t-a)}\,\Omega(t,a)
\qquad
\Om(t,a)\,\frac{\tht_\tau(t+a)}{\tht_\tau(t-a)}\,\frac{\tht_p(t+a)}{\tht_p(t-a)}=
\Om(-t,a).
\ee
We also have  $\Omega(t,\tau,p;a)=\Omega(t,p,\tau;a)$.

Given a one-variable phase function $\Omega(t,a)$, we define an $l$-variable \emph{phase function} by
\be\label{phasem}
\Omega(t_1,\dots,t_l,z_1,\dots,z_n, a_1,\dots,a_n)=\prod_{j=1}^l\prod_{m=1}^n\Omega(t_j-z_m, a_m)\prod_{1\leq i<j\leq l}\Omega(t_i-t_j, -2\eta).
\ee

\subsection {The elliptic action of the symmetric group}\label{actions}

Let $f=f(t_1,\dots,t_l)$ be a function. For complex numbers $\tau,\eta$ and a permutation $\sigma\in{\Bbb S}^l$, define the functions $[f]^{\tau,\eta}_{\sigma}$ via the action of the simple transpositions $(i,i+1)\in {\Bbb S}^l$, $i=1,\dots,\l-1$, given by      
\be
[f]^{\tau,\eta}_{(i,i+1)}(t_1,\dots, t_l)=
f(t_1,\dots,t_{i+1},t_i,\dots,t_l)\frac{\theta_\tau(t_i-t_{i+1}-2\eta)}{\theta_\tau(t_i-t_{i+1}+2\eta)}.
\ee

If for all $\sigma\in{\Bbb S}^l$, $[f]_{\sigma}=f$, we will say that the function is \emph{symmetric with respect to the elliptic action}.
 
\subsection {Spaces of elliptic weight functions}

Let $n=1$. Define \emph{a one-point elliptic weight functions} $w_l(t_1,\dots,t_l,\la,\tau,\eta;z,a)$ by
\be\label{eb}
w_l(t_1,\dots,t_l,\la,\tau,\eta;z,a)=
\prod_{i<j}\frac{\theta_\tau(t_i-t_j)}
{\theta_\tau(t_i-t_j+2\eta)}\prod_{j=1}^{l}\frac
{\theta_\tau(t_j-z-a+\lambda+2\eta l)}
{\theta_\tau(t_j-z-a)}.
\ee

Let $n,l$ be natural numbers. Set $\Zb^n_l=\{\bar{l}=(l_1,\dots,\l_n)\in\Z^n_{\ge 0}\,|\,\sum\limits_{i=1}^nl_i=l\}$ and for $m=0,1,\dots,n$, set $l^m=\sum\limits_{i=1}^ml_i$.

For $\bar{l}\in\Zb^n_l$, define \emph{an elliptic weight function} $w_{\bar{l}}(t_1,\dots,t_l,\la,\tau,\eta;z_1\dots,z_n,a_1\dots,a_n)$ by 

\bean\label{ w.f }
\lefteqn{w_{\bar{l}}(t,\la,\tau,\eta;z,a)=}
\\&&
\sum_{\sigma\in{\Bbb S}^l} \left[\prod_{i=1}^n \left(\frac{1}{l_i!} w_{l_i}(t_{l^{i-1}+1},\dots, t_{l^i}, \la-\sum_{j=1}^{i-1}2\eta\mu_j,\tau,\eta; z_i,a_i)\prod_{m=1}^{i-1}\prod_{s=l^{i-1}}^{l^i}
\frac{\th_\tau(t_m-z_s+a_s)}{\th_\tau(t_m-z_s-a_s)}\right)
\right]_\sigma,\notag
\eean
where $\mu_i=a_i/\eta-2l_i$, $i=1,\dots, n$.

For fixed $z,a\in\C^n$, $\la,\tau,\eta\in\C$, the space spanned over $\C$ by all elliptic weight functions $w_{\bar{l}}(t,\la,\tau,\eta;z,a)$, $\bar{l}\in\Zb_n^l$, is called \emph{the space of elliptic weight functions} and is denoted by $F_l(\la,\tau,\eta;z,a)$. Set $F_0(\la,\tau,\eta;z,a)=\C$ and $F(\la,\tau,\eta;z,a)=\bigoplus\limits_{l=0}^\infty F_l(\la,\tau,\eta;z,a)$.

Let $\h=\C h$ be a one-dimensional Lie algebra with generator $h$. For each
$\Lambda\in\C$ consider the $\h$-module $V_\Lambda=\oplus_{j=0}^\infty\C e_j$,
with $he_j=(\Lambda-2j)e_j$. We think of $V_\La$ as an evaluation Verma module over the quantum elliptic group $E_{\tau,\eta}(sl_2)$, see \cite{FV1}.

Let $V_\Lambda^*=\oplus_{j=0}^\infty\C e_j^*$ be the restricted dual of
the module $V_\La$. It is spanned by the basis
$(e_j^*)$ dual to the basis $(e_j)$. We let $\h$ act on $V_\Lambda^*$ by
$he_j^*=(\Lambda-2j)e_j^*$.  

Let $\om (\la,\tau,\eta;z,a):\, V_{\La_1}^*\T\dots\T V_{\La_n}^*\to F(\la,\tau,\eta;z,a)$ be a $\C$-linear map sending $e_{l_1}^*\T\dots\T e_{l_n}^*$ to $w_{\bar{l}}(t,\la,\tau,\eta;z,a)$. For generic values of parameters, the map $\om (\la,\tau,\eta;z,a)$ is an isomorphism of vector spaces, see \cite{FTV}.

\subsection{Admissible weight functions}\label{adm}

Let $\La=(\La_1,\dots,\La_n)\in\C^n$, $\bar{l}=(l_1,\dots,l_n)\in\Z^n_{\ge 0}$.
An $i$-th coordinate of $\bar{l}$ is called
\emph{$\La$-admissible} if either $\La_i\not\in\Z_{\ge 0}$
or $\La_i\in\Z_{\ge 0}$ and $l_i\le \La_i$.
The index $\bar{l}$ is called \emph{$\La$-admissible} if all its coordinates are $\La$-admissible. Denote $B_\La(\bar{l})\subset\{1,\dots,n\}$ the set of
all non-$\La$-admissible coordinates of $\bar{l}$.

For fixed $z,a\in\C^n$, $\la,\tau,\eta\in\C$, the space spanned over $\C$ by elliptic weight functions $w_{\bar{l}}(t,\la,\tau,\eta;z,a)$, $\bar{l}\in\Zb_n^l$, $\bar{l}$ is $\La$-admissible, is called \emph{the space of $\La$-admissible elliptic weight functions} and is denoted by $F^{\rm adm}_l(\la,\tau,\eta;z,a)$. Set $F^{\rm adm}_0(\la,\tau,\eta;z,a)=\C$ and $F^{\rm adm}(\la,\tau,\eta;z,a)=\bigoplus\limits_{l=0}^\infty F^{\rm adm}_l(\la,\tau,\eta;z,a)$.

For a non-negative integer $\La$, let $S_\La\in V_\La$ be the subspace given by $S_\La=\bigoplus\limits_{j=\La+1}^\infty\C e_j$. For a complex number $\La\not\in\Z_{\ge 0}$, let $S_\La$ be the trivial subspace. Let $L_\La=V_\La/S_\La$. For a non-negative integer $\La$, we think of $L_\La$ as a finite dimensional evaluation representation of the quantum elliptic group $E_{\tau,\eta}(sl_2)$, see \cite{FV1}. 

Define a map $\om ^{\rm adm}(\la,\tau,\eta;z,a):\, L_{\La_1}^*\T\dots\T L_{\La_n}^*\to F^{\rm adm}(\la,\tau,\eta;z,a)$ by the same formula as $\om (\la,\tau,\eta;z,a)$.

\section{The QKZB equation}\label{QKZB}

\subsection{The elliptic R-matrix}

Let $n=2$. The spaces of functions $F(\la,\tau,\eta;z_1,z_2,\La_1,\La_2)$ and $F(\la,\tau,\eta;z_2,z_1,\La_2,\La_1)$ coincide, see \cite{FTV}. 
The \emph{elliptic $R$-matrix} is an operator 
\newline
$R(\la,\tau,\eta; z_1,z_2,\La_1,\La_2)\in\End(V_{\La_1}\T V_{\La_2})$ dual to the transition matrix expressing the basis $\tilde w_{ij}=\om (\la,\tau,\eta;z_2,z_1,a_2,a_1)e_j^*\otimes e_i^*$ of the space $F(\la,\tau,\eta;z,\La)$ in terms of the basis $w_{ij}=\om (\la,\tau,\eta;z_1,z_2,a_1,a_2) e_i^*\otimes e_j^*$.
Namely, 
\be
R(\la,\tau,\eta;z_1,z_2,\La_1,\La_2)e_i\otimes e_j=
\sum_{kl}R_{ij}^{kl}e_k\otimes e_l, \qquad 
\tilde\omega_{kl}=\sum_{ij}R_{ij}^{kl}\omega_{ij}.
\ee
 
The elliptic $R$-matrix $R(\la,\tau,\eta;z_1,z_2,\La_1,\La_2)$ is characterized by an intertwining property with respect to the action of the elliptic quantum
group $E_{\tau,\eta}(sl_2)$ on tensor products of evaluation Verma modules, see \cite{FV1}.

The elliptic $R$-matrix has the following properties, see \cite{FTV}.
\begin{enumerate}
\item $R(\la,\tau,\eta;z_1,z_2,\La_1,\La_2)$ is a meromorphic function of $\la,z_1,z_2,\La_1,\La_2$.
\item
$R(\la,\tau,\eta;z_1,z_2,\La_1,\La_2)$ depends on $z_1,z_2$ only through the difference $z_1-z_2$. Accordingly, we write $R(\la,\tau,\eta;z_1-z_2,\La_1,\La_2)$ instead of
$R(\la,\tau,\eta;z_1,z_2,\La_1,\La_2)$.
\item The zero weight property: for all $\la,\tau,\eta,z,\La_1,\La_2$,
$[R(\la,\tau,\eta;z,\La_1,\La_2), h^{(1)}+h^{(2)}]=0$.
\item For any $\La_1,\La_2,\La_3$, the dynamical Yang--Baxter
equation
\be
R^{(12)}(\lambda-2\eta h^{(3)},\tau,\eta;z,\Lambda_1,\Lambda_2)
R^{(13)}(\la,\tau,\eta;z+w,\Lambda_1,\Lambda_3)
R^{(23)}(\lambda-2\eta h^{(1)},\tau,\eta;w,\Lambda_2,\Lambda_3)
\ee
\be
=
R^{(23)}(\la,\tau,\eta;w,\La_2,\La_3)
R^{(13)}(\la-2\eta h^{(2)},\tau,\eta;z+w,\La_1,\La_3)
R^{(12)}(\la,\tau,\eta; z,\La_1,\La_2)
\ee
holds in $\End(V_{\Lambda_1}\otimes V_{\Lambda_2}\otimes V_{\Lambda_3})$
for all $\la,\tau,\eta,z,w$.
\item For all $\la,\tau,\eta,z,\La_1$, $\La_2$, 
$R^{(12)}(\la,\tau,\eta;z,\La_1,\La_2)R^{(21)}(\la,\tau,\eta;-z,\La_2,\La_1)=\Id$.
This property is called \emph{the unitarity of $R$-matrix}.
\item The elliptic $R$-matrix $R(\la,\tau,\eta;z,\La_1,\La_2)$ preserves
$S_{\La_1}\T V_{\La_2}+V_{\La_1}\T S_{\La_2}$ for all $\la,z$. In particular, the elliptic $R$-matrix $R(\la,\tau,\eta;z,\La_1,\La_2)$ induces operators, still denoted by $R(\la,\tau,\eta;z,\La_1,\La_2)$, on the quotients $L_{\La_1}\T L_{\La_2}$. These operators obey the dynamical Yang--Baxter equation. 
\end{enumerate}

We use the following notation: if $X\in\End(V_i)$, then we denote by
$X^{(i)}\in\End(V_1\otimes\ldots\otimes V_n)$ the operator
$X$, acting non-trivially on
the $i$-th factor of a tensor product of vector spaces, and if
$X=\sum X_k\otimes Y_k\in\End(V_i\otimes V_j)$, then we set
$X^{(ij)}=\sum X_k^{(i)}Y_k^{(j)}$. If $X(\mu_1,\dots,\mu_n)$ is a function
with values in $\End(V_1\otimes\ldots\otimes V_n)$, then
$X(h^{(1)},\dots,h^{(n)})v=X(\mu_1,\dots,\mu_n)v$ if $h^{(i)}v=\mu_iv$,
for all $i=1,\dots,n$.

\subsection{The qKZB equations}
Fix complex parameters $\eta,\tau,p$. Fix also $n$ complex numbers
$\Lambda_1,\dots,\Lambda_n$.
Let $M=V_{\Lambda_1}\otimes\cdots\otimes V_{\Lambda_n}$ or $M=L_{\Lambda_1}\otimes\cdots\otimes L_{\Lambda_n}$. The kernel of
$h^{(1)}+\dots+h^{(n)}$ on $M$ is called the zero-weight subspace and is denoted
$M[0]$. The quantized Knizhnik-Zamolodchikov-Bernard (qKZB) equation is a system of difference equations for
a function $\Psi(\la,\tau,p,\eta;z_1,\dots,z_n)$ of $n$ complex variables
$z_1,\dots,z_n$ with values in the space of meromorphic functions of a complex variable $\la$ with values in $M[0]$, where $\tau,p,\eta$ are parameters of the equation.

The qKZB equation \cite{F} has the form
\be\label{qKZB}
\Psi(\la,\tau,p,\eta;z_1,\dots,z_j+p,\dots,z_n,a)= K_j(\la,\tau,p,\eta;z_1,\dots,z_n,a)\Psi(\la,\tau,p,\eta;z_1,\dots,z_n,a),
\ee
\be
K_j(\la,\tau,p,\eta;z,a)=R_{j,j-1}(z_j\!-\!z_{j-1}+p)\ldots R_{j,1}(z_j\!-\!z_{1}+p)\Gamma_jR_{j,n}(z_j\!-\!z_n)\cdots,R_{j,j+1}(z_j\!-\!z_{j+1}),
\ee
$j=1,\dots,n$. Here $R_{k,m}(z)$ is the elliptic $R$-matrix 
$R^{(k,m)}(\la-2\eta\sum h^{(s)},\tau,\eta;z,\La_k,\La_m)$, where the sum is taken over $s=1,\dots,m-1$, $s\neq k$,
acting on the $k$-th and $m$-th factors of the tensor product,
and $\Gamma_j$ is the linear difference operator such that
$\Gamma_j\Psi(\lambda)=\Psi(\lambda-2\eta\mu)$ if $h^{(j)}\Psi= \mu\Psi$.


\section{Analytic properties of hypergeometric integrals}\label{integrals}

\subsection{Hypergeometric integrals}

Fix natural numbers $n,l$ and $N$. Fix complex parameters $\tau$, $\eta$, $p$ such that $\Im\tau>0$, $\Im\eta<0$ and $\Im p>0$. Fix $\La=(\Lambda_1,\dots, \Lambda_n),\, z=(z_1,\dots, z_n)\in\C^n$, and set $a_j=\eta\La_j,\, j=1,\dots,n$ . Let $a=(a_1,\ldots,a_n)$, $t=(t_1,\ldots,t_l)$.
Let $\Om(t,\tau,p,\eta;z,a)$ be the $l$-variable phase function introduced in
\Ref{phase}. Fix $\bar{l},\bar{m}\in\Zb_n^l$ and let $w_{\bar{l}}(t,\la,\tau,\eta;z,a)$ and $w_{\bar{m}}(t,\mu,p,\eta;z,a)$
be the elliptic weight functions introduced in \Ref{ w.f }.
Let $\xi$ be an entire function of one variable
which is $4\eta N$-periodic, $\xi(\la+4\eta N)=\xi(\la)$.

Assume that for all $i=1,\ldots,n$,
\be\label{hor}
\Im a_i\gg\Im\tau\,,\quad\ \Im a_i\gg\Im p\,,\quad\ \text{and}\ \quad
-\Im\eta\gg\Im\tau\,,\quad\ -\Im\eta\gg\Im p\,.
\ee

Consider the integral
\bean\label{hyperg}
I_{\bar{l},\bar{m}}^\xi(\la,\mu,\tau,p,\eta;z,a)=\int_{[0,N]^l}\,X_{\bar{l},\bar{m}}^\xi(t,\la,\mu,\tau,p,\eta;z,a)\,dt,
\eean
where $dt=dt_1\wedge\ldots\wedge dt_l$ and 
\bea
\lefteqn{X_{\bar{l},\bar{m}}^\xi(t,\la,\mu,\tau,p,\eta;z,a)=}
\\ &&
\xi(p\la+\tau\mu-\Sum_{l=1}^n 2a_lz_l+4\eta\Sum_{j=1}^m t_j)\,
\Om(t,\tau,p,\eta;z,a)\,w_{\bar{l}}(t,\la,\tau,\eta;z,a)\,
w_{\bar{m}}(t,\mu,p,\eta;z,a).
\eea

The integral will be called {\it a hypergeometric integral}. We call the parameter $\Sum_{i=1}^n a_i/\eta-2l$ the \emph{weight of hypergeometric integral}.
Note that the integrand of the hypergeometric integral, $X_{\bar{l},\bar{m}}^\xi(t,\la,\mu,\tau,p,\eta;z,a)$, in \Ref{hyperg} is N-periodic.

In order to define the function $I_{\bar{l},\bar{m}}^\xi(\la,\mu,\tau,p,\eta;z,a)$
for other  values of parameters we use analytic continuation. The result of the analytic continuation can be represented as an integral of the integrand over a suitably deformed torus, which we denote by $T^l_N$. 

Namely, the poles of the integrand $X_{\bar{l},\bar{m}}^\xi(t,\la,\mu,\tau,p,\eta;z,a)$
are at most of first order and lie at the hyperplanes given by equations
\bean\label{poles1}
t_i-z_k=\pm(a_k+rp+s\tau)+m, \qquad t_i=t_j\pm(2\eta-rp-s\tau)+m,
\eean
where $1\leq i<j\leq l$, $k=1,\ldots,n$ and $r,s\,\in\Z_{\ge 0}$, $m\in\Z$.

We move the parameters $\tau,\,p,\,\eta$, $a_1,\ldots,a_n$, $z_1,\ldots,z_n$  and deform the integration torus accordingly so that the torus does not intersect the hyperplanes \Ref{poles1} at every moment of the deformation.
Then the analytic continuation of the function
$I_{\bar{l},\bar{m}}^\xi(\la,\mu,\tau,p,\eta;z,a)$ is given by the integral
\be\label{Hyperg.al}
I_{\bar{l},\bar{m}}^\xi(\la,\mu,\tau,p,\eta;z,a)\,=\,\int_{T^l_N}\,X_{\bar{l},\bar{m}}^\xi(t,\la,\mu,\tau,p,\eta;z,a)\,dt\,.
\ee

\begin{thm}\label{cite}(Theorem 29 in \cite{FTV}) The hypergeometric integral $I_{\bar{l},\bar{m}}^\xi(\tau,p,\eta;z,a)$ can be analytically continued as a holomorphic univalued function to the domain of the parameters $\la, \mu,\tau,p,\eta,a,z$ satisfying conditions \Ref{im}-\Ref{z's} below.$\;\Box$
\end{thm}

\begin{equation}\label{im}
\Im\tau>0\,,\qquad \Im p>0\,,\qquad \Im\eta<0\,.
\end{equation}
\begin{equation}\label{indep}
\text{The numbers $\tau$ and $p$ are linearly independent over $\Z$.}
\end{equation}
\begin{equation}\label{eta}
\{2\eta,\, 4\eta,\ldots,\,2m\eta\}\cap\{\Z+\tau \Z+p\Z\}\,=\,\varnothing\,.
\end{equation}
\bean\label{a's}
\{2a_k+2s\eta\}\cap\{\Z+\tau\Z+p\Z\}\,=\,\varnothing,\qquad k=1,\ldots,n\,,\qquad
s=1-m,\ldots,m-1\,.
\eean
\bean\label{z's}
\{z_l\pm a_l-z_k\pm a_k+2s\eta\}\cap\{\Z+\tau\Z+p\Z\}\,=\,\varnothing,\qquad
k,l=1,\ldots,n\,,\quad l\neq k\,,
\\
s=1-m,\ldots,m-1\,,
\notag
\eean
for arbitrary combination of signs.

\subsection{Analytic continuation}
In this Section we formulate a stronger version of Theorem~\ref{cite}.

We will often assume that the parameters $\tau, p,\eta, z_1,\dots,z_n, \La_1,\dots,\La_n$ satisfy the following conditions.

\bean\label{eta1}
\{2\eta s\, | \, s\in \Z_{>0}, s<
\max\{2, \Re\La_1,\dots,\Re\La_n\},\, s\le l\}\cap\{\Z+\tau \Z+p\Z\}\,=\,\varnothing\,,
\eean
\bean\label{a's1}
\{2a_k-2s\eta\,|\,s\in \Z_{>0},\,s<\Re\La_m,\, s<l\}\cap\{\Z+\tau \Z+p\Z\}\,=\,\varnothing,\qquad k=1,\ldots,n\,,
\eean
\bean\label{z's1}
\{z_k-z_m\pm(a_k+a_m)+2s\eta\} \cap\{\Z+\tau \Z+p\Z\}\,=\,\varnothing,\qquad  s=1-l,\ldots,l-1,
\\
k,m=1,\ldots,n\,,\qquad m\neq k.
\notag
\eean

Note that conditions \Ref{eta1}-\Ref{z's1} imply conditions \Ref{eta}-\Ref{z's}.

\begin{thm}\label{an.cont}
Let $\tau,p,\eta\in\C$ satisfy conditions \Ref{im}, \Ref{indep}. Let $z^0,a^0=\eta\La^0\in\C^n$ satisfy conditions \Ref{eta1}-\Ref{z's1}. 
Let  $\bar{l}, \bar{m}\in \Zb^n_l$. Assume that for all $i=1,\dots,n$, either the $i$-th coordinate of $\bar{l}$ or 
the $i$-th coordinate of $\bar{m}$ is $\La^0$-admissible, i.e.
$B(\bar{l})\bigcap B(\bar{m})=\varnothing$.
Then the hypergeometric integral $I_{\bar{l},\bar{m}}^\xi(\la,\mu,\tau,p,\eta;z,a)$ is holomorphic at $z^0,a^0$.
Moreover, there exists a contour of integration $T^l_N (z^0,a^0)\subset\C^l$ independent on $\bar{l},\bar{m}$, such
that for all 
$z, a$ in a small neighborhood of $z^0, a^0$, we
have
\be
I_{\bar{l},\bar{m}}^\xi(\la,\mu,\tau,p,\eta;z,a)\,=\,\int_{T^l_N(z^0,a^0)}\,X_{\bar{l},\bar{m}}^\xi(t,\la,\mu,\tau,p,\eta;z,a)\,dt\,.
\ee
\end{thm}

The proof of Theorem \ref{an.cont} is the same as the proof of Theorem 10 in \cite{MV}.

\subsection{A functorial property of hypergeometric integrals}

Let $a^0=\eta\La^0\in\C^n,\,\bar{l},\bar{m}\in \Zb^n_l$. Recall that $B_{\La^0}(\bar{l})\in\{1,\dots,n\}$ denotes the set of all non-$\La^0$-admissible coordinates of $\bar{l}$.
Let $B=B_{\La^0}(\bar{l})\cap B_{\La^0}(\bar{m})$.

Introduce $\tilde\La^0(B)=(\tilde\La_1^0,\dots,\tilde\La_n^0)\in\C^n$ by the rule: $\tilde\La^0_i=\La_i^0$ if $i\not\in B$ and $\tilde\La_i^0=-\La_i^0-2$ if
$i\in B$. Let $\tilde a^0(B)=\eta\tilde\La^0(B)$.

Let $\bar{l}^\prime(B)=(l_1^\prime,\dots,l_n^\prime)$, where
$\l_i^\prime=l_i$ if $i\not\in B$ and $l_i^\prime=l_i-\La_i-1$ if
$i\in B$. Similarly, let $\bar{m}^\prime(B)=(m_1^\prime,\dots,m_n^\prime)$, where
$m_i^\prime=m_i$ if $i\not\in B$ and $m_i^\prime=m_i-\La_i-1$ if 
$i\in B$.

We have $\l_1^\prime+\ldots +l_n^\prime=m_1^\prime+\ldots +m_n^\prime=
l-\sum\limits_{i\in B}\La_i-|B|$, where $|B|$ is the number of elements in the set $B$. We denote this number $l^\prime(B)$.

Let $B=\{b_1,b_2,\dots,b_{|B|}\}$. For a function $g(z,a)$, introduce a function
\be
\left.
(\res_{a=a^0}g)(z,a^0)=(\res_{a_{b_{|B|}}=a^0_{b_{|B|}}}\ldots\res_{a_{b_2}=a^0_{b_2}}\res_{a_{b_1}=a^0_{b_1}}g(z,a))\right|_{a_i=a_i^0,\; i\not\in B}.
\ee

\begin{thm}\label{functor}
Let $\tau,p,\eta\in\C$ satisfy conditions \Ref{im}, \Ref{indep}. Let $z^0,a^0=\eta\La^0\in\C^n$ satisfy conditions \Ref{eta1}-\Ref{z's1}. 
Let  $\bar{l}, \bar{m}\in \Zb^n_l$ and $B=B_{\La^0}(\bar{l})\cap B_{\La^0}(\bar{m})$. Then the hypergeometric integral $I_{\bar{l},\bar{m}}^\xi(\la,\mu,\tau,p,\eta;z,a)$ has a pole of order $|B|$ at $z^0,a^0$. Moreover,
\bean\label{res}
\lefteqn{(\res_{a=a^0} I_{\bar{l},\bar{m}}^\xi)(\la,\mu,\tau,p,\eta;z^0,a^0)=}
\\ &&
\left(\prod_{s=l^\prime(B)+1}^l\frac{\th_\tau(\la+2s\eta)\,\th_p(\mu+2s\eta)}{s}\right)\,
c_B(\tau,p,\eta,z^0,a^0)\, I_{\bar{l}^\prime,\bar{m}^\prime}^\xi(\la,\mu,\tau,p,\eta;z^0,\tilde a^0),\notag
\eean
where $c_B(\tau,p,\eta,z,a^0)$ is a nonzero holomorphic function at $z^0$ given below.
\end{thm}

Note that  $I_{\bar{l}^\prime,\bar{m}^\prime}^\xi(\la,\mu,\tau,p,\eta;z^0,\tilde a^0)$ is an $l^\prime(B)$-dimensional hypergeometric integral. We have 
$B_{\tilde\La^0}(\bar{l}^\prime)\cap B_{\tilde\La^0}(\bar{m}^\prime)=\varnothing$. By Theorem \ref{an.cont} the integral $I_{\bar{l}^\prime,\bar{m}^\prime}^\xi(\la,\mu,\tau,p,\eta;z^0,\tilde a^0)$ is well defined. 
Note also that the weights of hypergeometric integrals in the right hand side and in the left hand side in \Ref{res} are equal. 
 
Theorem \ref{functor} connects $l$- and $l^\prime(B)$-dimensional hypergeometric integrals of the same weight.

Now, we describe the function $c(\tau,p,\eta,z,a^0)$. For $j\in\{1,\dots,n\}$, $k\in\Z_{\ge 0}$, let
\be
N_{j,k}(\tau,p,\eta;z,a)=\prod_{s=0}^{k}\left(\prod_{i=0}^{j-1}\Om(z_i-z_j,a_i+a_j-2s\eta)\prod_{i=j+1}^{n}\Om(z_j-z_i,a_i+a_j-2s\eta)\right). 
\ee
For $k\in\Z_{\ge 0}$, let 
\be
y_k(\tau,p,\eta)=\frac{d_k(\eta)}{(k+1)!}\,
(\Om^\prime(-2\eta,-2\eta))^k\,
\prod_{i=0}^{k-2}\prod_{j=i+2}^k \Om(2(i-j)\eta,-2\eta)
\prod_{i=1}^{k-1}\Om(k\eta-2i\eta,k\eta),
\ee
where $\Om^\prime(t,a)$ denotes derivative with respect $t$ and
\be
\left.d_k(\eta)=\Om(t,k\eta)\Om(-t,k\eta)\right|_{t=k\eta}.
\ee
For $k\in\Z_{\ge 0}$, let 
\be
x_k(\tau,\eta)=\frac{1}{(k+1)!\,\th_\tau^\prime(0)}\prod_{i=0}^{k-1}\prod_{j=i+1}^k
\frac{\th_\tau(2(i-j)\eta)}{\th_\tau(2(i-j+1)\eta)}\prod_{i=1}^k\frac{1}{\th_\tau(2i\eta)}.
\ee
Then
\be
c_B(\tau,p,\eta,z,a^0)=\prod_{s\in B}  x_{\La^0_s}(\tau,\eta)\,x_{\La^0_s}(p,\eta)\,y_{\La^0_s}(\tau,p,\eta)\, N_{s,\La^0_s}(\tau,p,\eta;z,a^0). 
\ee

The proof of Theorem \ref{functor} is the same as the proof of Theorem 12 in \cite{MV}.

\section{Solutions of the qKZB equation}\label{solutions}
\subsection{Solutions of the qKZB equation with values in $V_{\Lambda_1}\otimes\ldots\otimes V_{\Lambda_n}$}
Fix complex numbers $\tau,\,p,\,\eta$, $\Lambda_1,\ldots,\Lambda_n$ and set
$a_i=\eta\Lambda_i$. Assume that the parameters $\tau,\,p,\,\eta$,
$a_1,\ldots,a_n$ satisfy conditions 
\Ref{im}\,-\,\Ref{indep}, \Ref{eta1}-\Ref{a's1} and 
$\Lambda_1+\ldots+\Lambda_n=2l$ for some positive integer $l$.

Fix a natural number $N$. Let $\xi$ be an entire function of one variable
which is $4\eta N$-periodic, $\xi(\la+4\eta N)=\xi(\la)$.

Let $V_\La=\oplus_{j=0}^\infty\C e_j$ and
$V=V_{\Lambda_1}\otimes\ldots\otimes V_{\Lambda_n}$.
For any $\bar{l}\in\Zb^l_n$, set
$e_{\bar{l}}=e_{l_1}\otimes\ldots\otimes e_{l_n}{{}\,\in\,V}$.

Introduce a $V[0]\otimes V[0]$-valued function $u^\xi$ by
\bean\label{u-sol}
u^\xi(\la,\mu,\tau,p,\eta;z,a)\,=\,e^{-\pi i{\mu\la\over2\eta}}\,\sum_{\bar{l},\bar{m}\in\Zb^l_n}\,
I^\xi_{\bar{l},\bar{m}}(\la, \mu, \tau,p,\eta;z,a)\, e_{\bar{l}}\T e_{\bar{m}}.\eean

Introduce a function
\be\label{alpha}
\al(\la,\eta)\,=\, \exp (-\pi i \la^2/4\eta)\,,
\ee
and operators $D(\la,\eta;a)\in\End(V[0])$

\bean\label{D}
D_{k}(\la,\eta;a)\,=\,
\frac{\al(\la\,-\,2\eta\,\sum_{m=1}^{k}h^{(m)},\eta)}{
\al(\la\,-\,2\eta\,\sum_{m=1}^{k-1}h^{(m)},\eta)}\,
\exp(\pi i\eta\Lambda_k(\sum_{m=1}^{k-1}\Lambda_m-\sum_{m=k+1}^{n}\Lambda_m)).
\eean

\begin{thm}\label{sol-cite}
Under the above conditions, for any $j=1,\ldots,n$, we have
\bea\label{eq1}
u^\xi(\ldots,z_j+p,\ldots)=K_j(\la,\tau,p;z_1,\ldots,z_n,a)
\otimes D^{-1}_j(\mu,\eta;a)\,\,u^\xi(\ldots,z_j,\ldots),
\eea
\bea\label{eq2}
u^\xi(\ldots,z_j+\tau,\ldots)=D_j^{-1}(\la,\eta;a)\otimes K_j(\mu,p,\tau; z_1,\ldots,z_n,a)
\,\,u^\xi(\ldots,z_j,\ldots),
\eea
and, if in addition the function $\xi$ is $2a_l$-periodic
for all $l=1,\ldots,n$, then
\bea\label{eq3}
u^\xi(\ldots,z_j+1,\ldots)=u^\xi(\ldots,z_j,\ldots).
\eea
Here $K_j(\la,\tau,p,\eta;z_1,\dots,z_n,a)$ is the qKZB operator defined in Section \ref{qKZB}, i.e.\ the operator of the qKZB equations with step $p$ and
defined in terms of elliptic $R$-matrices with modulus $\tau$.
\end{thm}

Theorem \ref{sol-cite} follows from Theorem 31 in \cite{FTV} and Theorem \ref{an.cont}.

Let $f:V[0]\to\C$ be a linear function.
Consider the functions
\bea\label{solu1}
\Psi_f^\xi(z,\la,\mu,\tau,p)\,=\,
(1\otimes f)\,(1\otimes D(\mu,p,\eta;z,a))\,u^\xi(\la,\mu,\tau,p;z,a)\,
\eea
and
\bea\label{solu2}
\Phi_f^\xi(z,\la,\mu,\tau,p)\,=\,
(f\otimes 1)\,(D(\la,\tau,\eta;z,a)\otimes 1)\,u^\xi(\la,\mu,\tau,p;z,a)\,,
\eea
where an $\End(V[0])$-valued function $D(\mu,p,\eta;z,a)$ is given by
\bean\label{De}
D(\mu,p,\eta;z,a)\,=\,\prod_{j=1}^n\,D_{j}(\mu,\eta;a)^{z_j/p}\,.
\eean

Theorem \ref{sol-cite} means that for a fixed $\mu$, the function $\Psi_f^\xi$ is a solution of the qKZB equations with modulus $\tau$ and step $p$, and for a fixed $\la$, the function $\Phi_f^\xi$ is a solution
of the qKZB equations with modulus $p$ and step $\tau$.

\subsection{Solutions of the qKZB equation with values in $L_{\Lambda_1}\otimes\ldots\otimes L_{\Lambda_n}$}\label{sol.qKZB} Let $L_\La=V_{\La}/S_{\la}$ as in Section \Ref{adm} and
$L=L_{\Lambda_1}\otimes\ldots\otimes L_{\Lambda_n}$.
Introduce a $L[0]\otimes L[0]$-valued function $u_{\rm adm}^\xi$ by
\bean\label{adm sol}
u_{\rm adm}^\xi(\la,\mu,\tau,p,\eta;z,a)\,=\,e^{-\pi i{\mu\la\over2\eta}}\,\sum\,I^\xi_{\bar{l},\bar{m}}(\la, \mu, \tau,p,\eta;z,a)\, e_{\bar{l}}\T e_{\bar{m}},
\eean
where the sum is over all $\La$-admissible indices $\bar{l},\bar{m}\in\Zb^l_n$.

\begin{thm}\label{th-sol}
Under the above conditions, for any $j=1,\ldots,n$, we have
\bea\label{adm-eq1}
u_{\rm adm}^\xi(\ldots,z_j+p,\ldots)=K_j(\la,\tau,p;z_1,\ldots,z_n,a)
\otimes D^{-1}_j(\mu,\eta;a)\,\,u_{\rm adm}(\ldots,z_j,\ldots),
\eea
\bea\label{adm-eq2}
u_{\rm adm}^\xi(\ldots,z_j+\tau,\ldots)=D_j^{-1}(\la,\eta;a)\otimes K_j(\mu,p,\tau;z_1,\ldots,z_n,a)
\,\,u_{\rm adm}^\xi(\ldots,z_j,\ldots),
\eea
and, if in addition the function $\xi$ is $2a_l$-periodic
for all $l=1,\ldots,n$, then
\bea\label{adm-eq3}
u^\xi_{\rm adm}(\ldots,z_j+1,\ldots)=u^\xi_{\rm adm}(\ldots,z_j,\ldots).
\eea
Here $D(\la,\eta;a)\in\End(L[0])$ is defined by \Ref{D} and $K_j(\la,\tau, p,\eta;z_1, \ldots,z_n,a)\in\End(L[0])$ is the qKZB operator defined in Section \ref{qKZB}, i.e.\ the operator of the qKZB equations with step $p$ and
defined in terms of elliptic $R$-matrices with modulus $\tau$.
\end{thm}

Theorem \ref{th-sol} follows from Theorem \ref{sol-cite}.

Let $f:L[0]\to\C$ be a linear function.
Consider the functions
\bea\label{adm-solu1}
\Psi_f^\xi(z,\la,\mu,\tau,p)\,=\,
(1\otimes f)\,(1\otimes D(\mu,p,\eta;z,a))\,u^\xi_{\rm adm}(\la,\mu,\tau,p;z,a)\,
\eea
and
\bea\label{adm-solu2}
\Phi_f^\xi(z,\la,\mu,\tau,p)\,=\,
(f\otimes 1)\,(D(\la,\tau,\eta;z,a)\otimes 1)\,u^\xi_{\rm adm}(\la,\mu,\tau,p;z,a)\,,
\eea
where an $\End(L[0])$-valued function $D(\mu,p,\eta;z,a)$ is given by \Ref{De}.

Theorem \ref{th-sol} means that for a fixed $\mu$, the function $\Psi_f^\xi$ is a solution of the qKZB equations with modulus $\tau$ and step $p$, and for a fixed $\la$, the function $\Phi_f^\xi$ is a solution
of the qKZB equations with modulus $p$ and step $\tau$.

\subsection{A functorial property of the hypergeometric solutions}
Let $\La^0\in\C^n$, $V_{\La^0}=V_{\La^0_1}\T \ldots V_{\La^0_n}$ and $L_{\La^0}=L_{\La^0_1}\T \ldots\L_{\La^0_n}$
For a subset $B\subset\{1,\dots,n\}$, introduce a $V[0]\T V[0]$-valued function $u_B^\xi$ by
\be\label{B sol}
u_{B}^\xi(\la,\mu,\tau,p,\eta;z,a)\,=\,e^{-\pi i{\mu\la\over2\eta}}\,\sum\,I^\xi_{\bar{l},\bar{m}}(\la, \mu, \tau,p,\eta;z,a)\, e_{\bar{l}}\T e_{\bar{m}},
\ee
where the sum is over all indices $\bar{l},\bar{m}\in\Zb^l_n$ such that $B_{\La^0}(\bar{l})= B_{\La^0}(\bar{m})=B$. 

Note that $u_B^\xi$ can be considered as a square submatrix of $u^\xi$ given by \Ref{u-sol}, in particular, we can identify $u_{B=\varnothing}^\xi(\la,\mu,\tau,p,\eta;z,a^0)= u_{\rm adm}^\xi(\la,\mu,\tau,p,\eta;z,a^0)$.

\begin{thm}\label{new-sol} 
Let $\tau,p,\eta\in\C$ satisfy conditions \Ref{im}, \Ref{indep}. Let $z^0,a^0=\eta\La^0\in\C^n$ satisfy conditions \Ref{eta1}-\Ref{z's1}. Then
\bea
\lefteqn{(\res_{a=a^0} u_B^\xi)(\la,\mu,\tau,p,\eta;z^0,a^0)=}
\\ &&
\left(\prod_{s=l^\prime(B)+1}^l\frac{\th_\tau(\la+2s\eta)\,\th_p(\mu+2s\eta)}{s}\right)\,
c_B(\tau,p,\eta,z^0,a^0)\, u_{\rm adm}^\xi(\la,\mu,\tau,p,\eta;z^0,\tilde a^0),
\eea
where $c_B(\tau,p,\eta,z,a^0)$ is a scalar nonzero holomorphic function at $z^0$ the same as in Theorem \ref{functor} and the parameter $\tilde a^0$ is also the same as in Theorem \ref{functor}.
\end{thm}

Theorem \ref{new-sol} follows from Theorem \ref{functor}.

\subsection{Monodromy of the solutions}
Set
\be
\Psi^\xi(\la,\mu,\tau,p,\eta;z,a)\,=\,(1\otimes D(\mu,p,\eta;z,a))\,
u^\xi_{\rm adm}(\la,\mu,\tau,p,\eta;z,a)\,,
\ee
where $u^\xi_{\rm adm}$ is given by \Ref{adm sol} and $D(\mu,p,\eta;z,a)$ is given by \Ref{De}.
According to Theorem \ref{adm sol}, the $L[0]\otimes L[0]$-valued function
$\Psi^\xi$ is a solution of the qKZB equation with respect to the first factor.

For $j=1,\ldots,n$, introduce a linear operator 
$B_j(\la,p,\eta;z,a)\in\End L[0]$ by
\be
B_j(\la,p,\eta;z,a)\,=
\,\exp(2\pi i\eta\Sum_{m,\,m\neq j}(z_m-z_j)\Lambda_m\Lambda_j/p)\,
D_j(\la,\eta;a).
\ee
Then the function
\be
\Psi^\xi_j(\la,\mu,\tau,p,\eta;z,a)\,=\,
(B_j(\la, p,\eta;z,a)\otimes 1)\,\Psi^\xi(\la,\mu,\tau,p,\eta;z_1,\dots,z_j+\tau,\dots,z_n,a)
\ee
is a new solution of the same equation, cf. Theorem 16 in \cite{FTV}.

The next Theorem describes a relation between the two solutions $\Psi^\xi_j$ and $\Psi^\xi$, and can be
considered as a description of the monodromy of the hypergeometric solutions
with respect to shifts of variables
$z_j$ by $\tau$.

\begin{thm}\label{mondr.shifts}
\bea\label{m.s}
\lefteqn{
\Psi^\xi_j(\la,\mu,\tau,p,\eta;z,a)\,=\,\bigl(1\,\otimes\,F_j(p,\eta;z,a)\,D(\mu,p,\eta;\ldots,z_j+\tau,\ldots,a)\times}
\\&&
\times
K_j(\mu,p,\tau,\eta;z,a)\,D^{-1}(\mu,p,\eta;\ldots,z_j,\ldots,a)\bigr)\,
\Psi^\xi(\la,\mu,\tau,p,\eta;z,a)
\eea
where
$F_j(p,\eta;z,a)=e^{2\pi i\eta\sum_{l,\,l\neq j}(z_l-z_j)\Lambda_l\Lambda_j/p}$
and
$K_j(\mu,p,\tau,\eta;z,a)$ is the $j$-th operator of the qKZB equations with step $\tau$ and modulus $p$.
\end{thm}

Theorem \ref{mondr.shifts} follows from Theorem 38 in \cite{FTV} and Theorem \ref{th-sol}. 

The monodromy of solution $u_{\rm adm}^\xi$ with respect to shifts  $z_i$ by $1$ is described by Proposition 39 in \cite{FTV} whose proof can be applied to our situation.  

The monodromy of solutions $u_{\rm adm}^\xi$ with respect to permutation of variables $z_i$ and $z_j$ is described by Theorem 36 in \cite{FTV} whose proof can be applied to our situation.   

The theta function properties of solutions $u_{\rm adm}^\xi$ are described by Theorem 33 in \cite{FTV} whose proof can be applied to our situation.

\end {document}